\documentclass[11pt]{amsart}
\usepackage{amssymb}
\newfont{\ams}{msbm10 at 12pt}
\newfont{\amsi}{msbm8}
\newcommand{\nc}{\newcommand}
\nc{\C}{{\mathcal C}}
\nc{\CB}{{\mathcal B}}
\nc{\supp}{\operatorname{supp}}
\nc{\M}{{\mathcal M}}
\nc{\LL}{{\mathbf L}}

\nc{\D}{{\mathbb D}}

\nc{\St}{\operatorname{St}^{\bullet}}
\nc{\A}{{\frak A}}

\nc{\1}{{\mathbf 1}}
\nc{\CC}{{\mathbb C}}
\nc{\R}{{\mathbb R}}
\renewcommand{\k}{{\mathbb  k}}
\nc{\Q}{{\mathbb Q}}
\nc{\U}{{\mathbf U}}
\nc{\B}{{\mathbf B}}
\nc{\N}{{\mathbb N}}
\nc{\Z}{{\mathbb Z}}
\nc{\Rhom}{\operatorname{RHom}\bul}
\nc{\Ad}{\operatorname{Ad}}
\nc{\Res}{\operatorname{Res}}
\nc{\gr}{\operatorname{gr}}
\nc{\tr}{\operatorname{tr}}
\nc{\End}{\operatorname{End}}

\renewcommand{\c}{{\mathfrak c}}

\nc{\g}{{\mathfrak g}}
\nc{\hatg}{\hat{\frak g}}

\renewcommand{\b}{{\mathfrak b}}
\nc{\sem}{{$\S_{\g}^{\g_{>0}}}}
\nc{\gl}{{\frak g\frak l}}
\nc{\n}{{\frak n}}
\nc{\si}{{\frac\infty 2}}
\nc{\p}{{\frak p}}
\nc{\h}{{\frak h}}

\nc{\Ind}{\operatorname{Ind}}
\nc{\ch}{\operatorname{ch}}
\nc{\Coind}{\operatorname{Coind}}
\nc{\opp}{{\operatorname{opp}}}
\nc{\Ker}{\operatorname{Ker}}
\nc{\im}{\operatorname{Im}}
\nc{\Coker}{\operatorname{Coker}}
\nc{\dirlim}{\underset{\rightarrow}{\operatorname{lim}}}
\nc{\invlim}{\underset{\leftarrow}{\operatorname{lim}}}
\nc{\Sem}{{{\mathbf S}_{\g}^{\g_{>0}}}}
\nc{\CN}{{\mathcal N}}
\nc{\Ext}{\operatorname{Ext}^{\bullet}}
\nc{\ext}{\operatorname{Ext}}
\nc{\tilW}{\til{W}}

\nc{\lth}{\ell t}
\nc{\BB}{\mathcal{B}}
\nc{\Tor}{\operatorname{Tor}_{\bullet}}
\nc{\tor}{\operatorname{Tor}}
\nc{\Tors}{\operatorname{Tor}_{\frac \infty 2+\bullet}}
\nc{\Exts}{\operatorname{Ext}^{\frac \infty 2+\bullet}}
\nc{\Hom}{\operatorname{Hom}^{\bullet}}
\nc{\ad}{\operatorname{ad}}
\renewcommand{\hom}{\operatorname{Hom}}

\renewcommand{\mod}{\operatorname{mod}}
\nc{\Mod}{\operatorname{Mod}}

\nc{\Barb}{\operatorname{Bar}^{\bullet}}

\nc{\upX}{X^{\uparrow}}
\nc{\upcD}{{\mathcal D}^{\uparrow}}
\nc{\upD}{D^{\uparrow}}
\nc{\dX}{X^{\downarrow}}
\nc{\dcD}{{\mathcal D}^{\downarrow}}
\nc{\dD}{D^{\downarrow}}
\nc{\upC}{{\mathcal C}^{\uparrow}}
\nc{\dC}{{\mathcal C}^{\downarrow}}
\nc{\underA}{\underline{A}}
\nc{\underC}{\underline{\CC}}
\nc{\underB}{\underline{B}}
\nc{\underk}{\underline{\k}}
\nc{\Db}{D^{\bullet}}
\nc{\ten}{{\otimes}}
\nc{\tenl}{\overset{\operatorname{L}}\ten}
\nc{\map}{\longrightarrow}

\nc{\bs}{\bigskip\\}
\nc{\ms}{\smallskip\\}

\nc{\tilbar}{\widetilde{\operatorname{Bar}}}
\nc{\tilBarb}{\widetilde{\operatorname{Bar}}^{\bullet}}
\nc{\overr}{\overline{R}}
\nc{\overI}{\overline{I}}
\nc{\overX}{\overline{X}}
\nc{\overh}{\overline{h}}
\nc{\overY}{\overline{Y}}
\nc{\overW}{\overline{W}}
\nc{\linbar}{\overline{\operatorname{Bar}}}

\nc{\til}{\widetilde}
\nc{\oppA}{A^{\sharp}}
\nc{\Lemma}{{\bf Lemma:\ }}
\nc{\Theorem}{{\bf Theorem:}\ }
\nc{\Cor}{{\bf Corollary:}\ }
\nc{\Def}{{\bf Definition:}\ }
\nc{\Prop}{{\bf Proposition:}\ }
\nc{\Con}{{\bf Conjecture:}\ }
\nc{\Rem}{{\bf Remark:}\ }
\nc{\dok}{{\bf Proof.}\ }

\nc{\SInd}{\operatorname{S-Ind}}
\nc{\SCoind}{\operatorname{S-Coind}}
\nc{\bul}{^{\bullet}}
\nc{\stand}{C^{\frac\infty2+\bullet}}
\nc{\ssn}{\subsection{}}
\nc{\sssn}{\subsubsection{}}

\nc{\hgt}{\operatorname{ht}}
\nc{\sqbinom}{\fracwithdelims[][0pt]}

\address{Independent University of Moscow, Pervomajskaya st. 16-18,
Moscow 105037, Russia}
\email{hippie@mccme.ru}
\author{Sergey Arkhipov}
\title{Semiinfinite cohomology of Tate Lie algebras}
\date{}

\begin{document}
\maketitle
\section{Introduction.}
This note is a natural extension of the final part  in \cite{Ar1}
 where a natural homological construction for graded Lie algebra
semiinfinite cohomology and a natural explanation for the phenomenon
of the critical cocycle was found.  Here we propose a variant of the
construction that works in the Tate Lie algebra case. Note that the
standard complex for the computation of the Tate Lie algebra
semiinfinite cohomology was written down by Beilinson and Drinfeld in
\cite{BD1} and probably by some physicists.  Still the construction
was rather an indirect one and was not formulated in terms similar to
the ones in \cite{Ar1}. In this note we spell out the construction of
the complex in terms of some kind of quadratic-linear Koszul duality
rather than Conformal Field Theory.

Let us say a few words about the contents of the note. In the second
section we recall the notion of a differential graded Lie algebra with a
curvature. We show that the standard Chevalley complex for computation of
the Lie algebra cohomology can be modified in such a way that it still
exists in this exotic case.

In the third section we recall a construction of an analogue of the
Chevalley cohomological complex for a {\em left} module over a Lie
algebroid $A$ over a commutative algebra $R$. Then we combine the
setup from the previous section with the described one and obtain a
picture consisting of a graded supercommutative algebra with a
derivation, a graded (super) Lie algebroid over the algebra carrying
an extension of the derivation of the basic algebra and finaly an
analogue of curvature that ``corrects'' the fact that both derivations
do not satisfy the constraint $d^2=0$. We construct an analogue of the
cohomological Chevalley complex with coefficients in a left CDG-module
over the differential graded Lie algebroid with curvature (CDG Lie
algebroid).

In the fourth section we show that the standard complex for
semiinfinite cohomology in the Tate Lie algebra case is a particular
example of the described situation.  Namely for a Tate Lie algebra
$\g$ with a compact Lie subalgebra $\b$ we consider the graded
supercommutative algebra $\Lambda\bul((\g/\b)^*)$ and a graded Lie
algebroid $\b\ten\Lambda\bul((\g/\b)^*)$ over it. We show that the
components of the Lie bracket in $\g$ provide the derivations and the
curvature.

Finally  we show that for a discrete module
$M^\sharp$ over the extension of $\g$ with the help of the cricical
cocycle of $\g$ the space $M^\sharp\ten \Lambda\bul(\g/\b)$ carries a
structure of a left CDG module over the above CDG Lie algebroid and
that for this CDG-module the standard Chevalley complex
from the third section coincides with the semiinfinite complex of
$\g$ with coefficients in $M^\sharp$.

In order to simplify the exposition we never use the language of
derived categories in the note but rather work with concrete
complexes.

\section{Toy example.}
The material of this section is based on a partly unpublished
construction of A.Polishchuk and L.Positselsky (still see \cite{P}).

We begin with the case of a differential graded Lie superalgebra with
a curvature.  The standard complex appearing in the case seems to be
more understandable.

Let $A=\oplus A_k$ be a graded Lie algebra and $d:\ A_k\map A_{k+1}$
be a derivation  of $A$ of  order $1$. Note that at this point we do
not put the constraint $d^2=0$. Instead we require an additional part
of the data --- an element $h\in A_2$ such that $d^2(a)=[h,a]$ for any
$a\in A$ and $d(h)=0$.

\ssn
\Def
The data $(A,d,h)$ described above are called the differential graded
Lie superalgebra with curvature or, for short, the CDG Lie algebra.

By definition a left (resp. a right) CDG-module over a CDG Lie algebra
$A$ is a graded left (resp. right) module $M=\oplus M_k$ over the Lie
algebra $A$ with the differential $d:\ M_k\map M_{k+1}$ satisfying the
Leibnitz rule such that $d^2(m)=h\cdot m$ (resp.  $d^2(m)=m\cdot h$)
for any $m\in M$.

Denote the category of left (resp. right) CDG  $A$-modules by
$CDG\operatorname{-}A\operatorname{-}\mod$
(resp.
$CDG\operatorname{-}\mod\operatorname{-}A$.

\subsection{Construction of the standard complex.} \label{cdglie}
For $M\bul\in CDG\operatorname{-}A\operatorname{-}\mod$ consider the
bigraded vector space $C^{\bullet\bullet}(A,M\bul)$ as follows:
$ C^{\bullet\bullet}(A,M\bul)= \hom(\Lambda\bul(A), M\bul) $,
here the first grading comes from the number of wedges in the exterior
product and
$$C^{\bullet k}(A,M\bul)=\underset{-p+q=k}\bigoplus
\hom((\Lambda\bul(A))_p,M^q).$$
Consider the two differentials on the bigraded vector space.  The
first one  of the grading $(1,0)$ is the usual Chevalley differential:
\begin{gather*}
(d_1f)(a_1\wedge\ldots\wedge a_m)=
\sum_i(-1)^if(a_1\wedge\ldots\wedge da_i\wedge\ldots\wedge a_m)
+d_{M\bul}
f(a_1\wedge\ldots\wedge a_m)\\
+\sum_i(-1)^ia_i^{M\bul}
f(a_1\wedge\ldots\wedge \hat a_i\wedge\ldots\wedge a_m)+
\sum_{i<j}(-1)^{i+j}f([a_i,a_j]\wedge a_1\wedge\ldots\wedge \hat a_i\wedge
\ldots\wedge \hat a_j\wedge\ldots\wedge a_m).
\end{gather*}
The second differential  of the grading $(-1,2)$ is given by the formula
that uses the curvature $h$:
$$
(d_2f)(a_1\wedge\ldots\wedge a_{m-1})=
f(h\wedge a_1\wedge\ldots\wedge a_{m-1}).
$$
Consider the total grading on the bigraded space.

\sssn                     \label{liedif}
\Lemma
The differential $d=d_1+d_2$ satisfies  $d^2=0$.

\dok
The corresponding calculation repeats  the one of Polishchuk and Positselsky.
\qed

\section{Standard complex for a CDG Lie algebroid.}
\subsection{Chevalley complex of a Lie algebroid.}
Let $R$ be a commutative algebra, and let $A$ be a Lie algebroid over
$R$, i.e. $A$ is a $R$-module carrying a Lie algebra structure over
the base field and acting on $R$ by derivations such that
$[a,rb]=a(r)b+r[a,b]$ for any $a,b\in A$ and $r\in R$.

By definition a $A$-module  is a $R$-module $M$ with the Lie
action of $A$ satisfying the constraint $a\cdot(rm)=r(a\cdot
m)+(a(r))\cdot m$ for any $a\in A$, $r\in R$ and $m\in M$.

\sssn
For a {\em left} $A$-module $M$
consider the graded vector space $C\bul(A,M)$ as follows:
$$C\bul(A,M)=\underset{k}\bigoplus C^k(A,M),\
C^k(A,M)=\hom_R(\Lambda_R^k(A),M)$$
We endow the graded vector space with the differential as follows:
\label{algebroid}
\begin{gather*}
(df)(a_1\wedge\ldots\wedge a_m)=
+\sum_i(-1)^ia_i^{M}
f(a_1\wedge\ldots\wedge \hat a_i\wedge\ldots\wedge a_m)\\
+\sum_{i<j}(-1)^{i+j}f([a_i,a_j]\wedge a_1\wedge\ldots\wedge \hat a_i\wedge
\ldots\wedge \hat a_j\wedge\ldots\wedge a_m).
\end{gather*}
\label{algebrcor}
\Lemma
\begin{itemize}
\item[(i)] The differential in the complex  is well defined.
\item[(ii)] The differential satisfies $d^2=0$.
\end{itemize}
\dok
(i) Let us perform a calculation showing that the differential
$d:\ C^k\map C^{k+1}$ is
well defined for $k=1$, the general case is quite similar.
Note that
\begin{gather*}
df(ra_1\wedge a_2)=
ra_1 f(a_2)
-a_2f(ra_1)-
f([ra_1,a_2])=
ra_1 f(a_2)
-a_2rf(a_1)+
-rf([a_1,a_2])+a_2(r)f(a_1)\\
=
ra_1 f(a_2)
-ra_2f(a_1)+
-rf([a_1,a_2])=
df(a_1\wedge ra_2).
\end{gather*}
(ii) The calculation does not differ from the usual Lie algebra case.
\qed

\sssn
\Rem Note that there is {\em no such standard complex} for a {\em right}
$A$-module $M$. \label{remark}

The constructed complex is called the (cohomological) Chevalley complex
of the Lie algebroid $A$ with coefficients in the left $A$-module $M$.

\subsection{CDG-algebroids.}
Now let $R=\oplus_kR_k$ be a graded supercommutative algebra and let
$A=\oplus A_k$ be a graded super Lie algebroid over $R$. Suppose also
there are a super derivation $d_R$ on $R$ of degree $1$ and a super
derivation $d_A$ of the Lie superalgebra $A$ also of degree $1$
satisfying Leibnitz rule with respect one to another. Again we put the
constraint  $d^2=0$ on neither of the derivations. Instead of that we
fix the choice of an element $h\in A_2$ such that $d_A^2(a)=[h,a]$ and
$d_R^2(r)=h(r)$.

\sssn
\Def
The data $(A,R,d_A,d_R,h)$ are called the differential graded Lie algebroid
with curvature or, for short, a CDG Lie algebroid.

The notion of a left (resp. right) CDG-module over a CDG algebroid
is a natural  combination of the previous definitions and we do not
spell it out explicitly.
The category of left (resp. right) CDG-modules ovea a CDG-algebroid
$A$ is denoted by
$CDG\operatorname{-}A\operatorname{-}\mod$
(resp.
$CDG\operatorname{-}\mod\operatorname{-}A$.

\subsection{Standard complex for a CDG Lie algebroid.}
Now we sort of put togeather definitions of the standard complexes
given in \ref{cdglie} and \ref{algebroid}.
For
$M\bul\in
CDG\operatorname{-}A\operatorname{-}\mod$
consider the bigraded vector space
$C^{\bullet\bullet}(A,M\bul)$ as follows:
$
C^{\bullet\bullet}(A,M\bul)=
\hom_R(\Lambda_R\bul(A), M\bul)
$,
here the first grading comes from the number of wedges in the exterior
product and
$$C^{\bullet k}(A,M\bul)=\hom_R^k((\Lambda_R\bul(A)),M\bul)$$
in the graded $\hom$ sense.
Consider the two differentials on the bigraded vector space.
The first one  of the grading $(1,0)$ is the usual Chevalley differential
like in \ref{algebroid}:
\begin{gather*}
(d_1f)(a_1\wedge\ldots\wedge a_m)=
\sum_i(-1)^if(a_1\wedge\ldots\wedge da_i\wedge\ldots\wedge a_m)
+d_{M\bul}
f(a_1\wedge\ldots\wedge a_m)\\
+\sum_i(-1)^ia_i^{M\bul}
f(a_1\wedge\ldots\wedge \hat a_i\wedge\ldots\wedge a_m)
+\sum_{i<j}(-1)^{i+j}f([a_i,a_j]\wedge a_1\wedge\ldots\wedge \hat a_i\wedge
\ldots\wedge \hat a_j\wedge\ldots\wedge a_m).
\end{gather*}
The second differential  of the grading $(-1,2)$ is again given by the formula
that uses the curvature $h$:
$$
(d_2f)(a_1\wedge\ldots\wedge a_{m-1})=
f(h\wedge a_1\wedge\ldots\wedge a_{m-1}).
$$
Consider the total grading on the bigraded space.

\sssn
\Lemma
\begin{itemize}
\item[(i)]
The differential $d_1$ is well defined and its square equals zero.
\item[(ii)]
The differential $d_2$ is well defined and its square equals zero.
\item[(iii)]
The differential $d=d_1+d_2$ satisfies  $d^2=0$.
\end{itemize}

\dok
(i) Follows from Lemma~\ref{algebrcor}.
(ii) Follows from the obvious  fact that $h\wedge h=0\in\Lambda^2_R(A)$.
(iii) Repeats the proof of Lemma~\ref{liedif}.
\qed

The obtained complex is called the (cohomological) Chevalley complex
of the CDG Lie algebroid $A$ with coefficients in the left CDG
$A$-module $M$.                       \label{maincomp}

\section{Semiinfiite cohomology via CDG Lie algebroids.}
In this section we show that the standard complex for the computation
of semiinfinite cohomology of a discrete module over a Tate Lie
algebra coincides with  the Chevalley  complex of the
form~\ref{maincomp} for a certain CDG Lie algebroid and a certain left
module over it.

\subsection{Tate Lie algebras.}
Recall that a {\em Tate} space is a complete topological vector space
having a base of neighbourhoods of $0$ consisting of commensurable vector
spaces(i.e., $\dim U_1/(U_1\cap U_2)<\infty$ for any $U_1$ and $U_2$
from this base).

Recall also that a {\em c-lattice} in a Tate space $V$ is an open
bounded subspace and a {\em d-lattice} $L\subset V$ is  a discrete
subspace such that there exists a c-lattice $P$ with $L+P=V$. Note that
the quotient of a Tate space by a c-lattice (resp. by a d-lattice)
is discrete (resp. compact) in its natural topology. It is known that
there is a natural duality on the category of Tate spaces and that
$V^{**}\til{\map}V$ for any Tate space V.

Recall also that by definition a {\em Tate Lie algebra} is a Tate
vector space equipped with a Lie algebra structure continuous in the
Tate topology.

\subsection{Construction pf the CDG Lie algebroid.}
Let $\g$ be a Tate Lie algebra with a subalgebra $\b\subset \g$ that
is a c-lattice In particular the space $\c:=\g/\b$ is discrete and its
dual space is compact. Choose a section of the projection $\g\map\c$.
Thus we fix a noncanonical decomposition $\g=\b\oplus\c$ and the Lie
algebra structure on $\g$ is provided by the following collection of
maps:
$$
\mu_\b:\ \b\wedge\b\map\b,\ \mu_1:\ \b\ten\c\map\c,\
\mu_\c:\ \c\wedge\c\map\c,\ \mu_2:\ \b\ten\c\map\b,\
h^*:\ \c\wedge\c\map\b.
$$
The existence of $h^*$ means exactly that $\c$ is not a Lie subalgebra
in $\g$.  The construction below is parallel to the one in \cite{Ar1}.

\sssn
First consider the graded supercommutative algebra $\Lambda\bul(\c^*)$
or rather its Tate completion denoted in the same way. So as a vector
spase it  is Tate dual to the discrete coalgebra $\Lambda\bul(\c)$.
Consider the derivation $d_R$ on $\Lambda\bul(\c^*)$ of degree one
generated by the map dual to $\mu_\c$ and extended to the whole
algebra by Leibnitz rule and by continuocity.

Next consider the graded Lie algebra $\b\ten\Lambda\bul(\c^*)$ where
the tensor product is understood in the completed sense so that the
whole space is Tate dual to the discrete space
$\b^*\ten\Lambda\bul(\c)$. The commutator map in the above Lie algebra
is generated by the one on $\b$ and by the action of $\b$ on the space
$\c^*=(\g/\b)^*$.

Denote the graded Lie algebra (resp. the graded supercommutative
algebra) above by $A$ (resp. by $R$). Evidently $A$ is a $R$-module,
moreover there is a natural adjoint action of $A$ on $R$.  \vskip 1mm
\noindent
\Lemma
$A$ is a graded Lie algebroid over the graded
supercommutative algebra $R$.  \qed

\subsubsection{Construction of the CDG Lie algebroid structure on $(A,R)$.}
We extend the derivation of $R$ constructed above to the derivation of
$A$. Namely consider the map $\b\map\b\ten\c^*$ dual to $\mu_2$.
Denote the map by $d_\b$. Extend the sum of $d_\b$ and $d_R$ to the
whole Lie algebra $A$ by Leibnitz rule and by continuocity. Denote the
obtained derivation by $d_A$. Finally consider the element $h\in
A_2=\b\ten\Lambda^2(\c^*)$ corresponding to the component $h^*$ of the
bracket in $\g$.

\sssn
\Prop
The data $(A,R,d_A,d_R,h)$ form a CDG Lie algebroid.
\qed

\subsection{Construction of the standard semiinfinite complex.}
First note that for any compact  $\g$-module $M$ the (completed)
tensor product $M\ten\Lambda\bul(\c^*)$ has a natural structure of a
{\em left} CDG-module over the CDG algebroid $(A,\ldots)$ constructed
in the previous subsection.

\subsubsection{Construction of the {\bf right} CDG-module.}
Now fix a {\em discrete} $\g$-module $M$ and consider the graded space
$M\ten\Lambda\bul(\c)$.
\vskip 1mm
\noindent
\Lemma
The graded space
$M\ten\Lambda\bul(\c)$ has a natural structure of the {\em right}
CDG-module over the CDG Lie algebroid $(A,\ldots)$.\qed

\sssn
Here we come to the crucial point explaining the phenomenon of the critical cocycle
in the semiinfinite cohomology of Tate Lie algebras. What we would like to do is to consider
the standard complex of the CDG Lie algebroid $(A,R,\ldots)$ with coefficients
in $M\ten\Lambda\bul(\c)$. Yet as noted in \ref{remark}
there is no naive way to do it. Somehow we have to make
$M\ten\Lambda\bul(\c)$ into a {\em left} CDG-module over our CDG Lie algebroid.

Happily the answer how to do this is known at least in the graded Lie
algebra case. The idea works for Tate Lie algebra semiinfinite cohomology
as well.

\subsubsection{Critical cocycle of $\g$.}
It is known that the Lie algebra of continuous endomorphisms of  $\g$
(and of any other Tate vector space as well) has a remarkable class in
$H^2$ and the choice of the  decomposition $\g=\b\oplus\c$ fixes its
representative in the space of cocycles called {\em the critical
2-cocycle} and denoted by $\omega_0$. The adjoint action of $\g$ on
itself provides the inverse image of the class denoted by
$\omega_0^\g$ and called the critical 2-cocycle of $\g$.

Denote by $\g^\sharp$ the central extension of $\g$ with the help of
this cocycle. Note that $\g^\sharp$ is a Tate Lie algebra with a fixed
vector space decomposition $\g=\b\oplus\c\oplus\CC K$, where $K$ denotes the
central element.

\sssn
Consider the CDG Lie algebroid $(\til A^\sharp,  R^\sharp,\ldots)$
obtained from the pair $\b\oplus \CC K\subset\g^\sharp$ by our main
construction.  Note that the algebra $R^\sharp\til\map
R=\Lambda\bul(\c^*)$ and that the  element $K$ is central in $\til A^\sharp$.
Denote by $A^\sharp$ the quotient of $\til A^\sharp$ by $K-1$.
\vskip 1mm
\noindent
\Prop
The CDG Lie algebroid
$(A^\sharp,R^\sharp,d_{A^\sharp},d_{R^\sharp},h^\sharp)$
is isomorphic to the
the CDG Lie algebroid
$(A,R,d_{A},d_{R},h)^{\operatorname{opp}}$.
\qed
\vskip 1mm
\noindent
\Rem
In particular for any discrete $\g^\sharp$-module $M^\sharp$ such that
the central element $K$ acts on it by $1$ the {\em right}
$(A^\sharp,\ldots)$ CDG-module $M\sharp\ten\Lambda\bul(\c)$
becomes a {\em left} CDG-module over the CDG Lie algebroid
$(A,\ldots)$.

\sssn
\Def
For a discrete $\g^\sharp$-module $M^\sharp$ such that
the central element $K$ acts on it by $1$
the complex
$$C\bul(A,M^\sharp\ten\Lambda(\c))=\hom_R(\Lambda_R\bul(A),
M^\sharp\ten\Lambda\bul(\c))
$$
is called the standard semiinfinite complex with coefficients in $M^\sharp$
and is denoted by
$C^{\si+\bullet}(\g,M^\sharp)$.

\sssn
\Theorem
\begin{itemize}
\item[(i)]
For a discrete $\g^\sharp$-module $M^\sharp$
as above the graded vector space
$C^{\si+\bullet}(\g,M^\sharp)$  is isomorphic to the graded vector space
$\Lambda\bul(\b^*)\ten\Lambda\bul(\g/\b)\ten M^\sharp$.
\item[(ii)]
The cohomology of the complex
$C^{\si+\bullet}(\g,M^\sharp)$ coincides with the semiinfinite cohomology
of the Tate Lie algebra $\g$ with coefficients in the discrete module
$M^\sharp$ defined in \cite{BD1}.
\qed
\end{itemize}


\begin{thebibliography}{99}
\bibitem[Ar]{Ar1}
S.M.Arkhipov. {\it
Semiinfinite cohomology of associative algebras and bar duality.}
International Math. Research Notices No. 17 (1997), 833-863.
\bibitem[BD]{BD1} A.Beilinson, V.Drinfeld. {\it
Quantizatiion of the Hitchin's integrable system and Hecke
eigensheaves.} Preprint, (1999), 1-384.
\bibitem[P]{P} L.Positselsky. Private communications.
\end{thebibliography}
\end{document}